\documentclass[12pt]{amsart}

\usepackage{amssymb}
\usepackage{amsmath,amsfonts}
\usepackage[latin1]{inputenc}

\usepackage{graphicx}
\usepackage{caption}

\makeatletter
\@namedef{subjclassname@2010}{%
  \textup{2010} Mathematics Subject Classification}
\makeatother

\usepackage[T1]{fontenc}
\newtheorem{thm}{Theorem}[section]

\newtheorem{Prop}{Proposition}[section]

\newcommand{\N}{\mathbb{N}}

\begin{document}

%%%%% To ease editing, for IMPAN journals add:

\baselineskip=17pt

%%%%%%%%%%%%%%%%

\title[The dynamic of q-deformed logistic maps]{Revisiting the dynamic of q-deformed logistic maps}

\author[Jose S. C\'anovas$^{1}$ and H. E. Rezgui$^{2}$]{Jose S. C\'anovas$^{1}$ and Houssem Eddine Rezgui$^{2}$}
\address{Jose S. C\'anovas$^{\text{1}}$, Department de Matemática Aplicada y estadística, Universidad Politécnica de Cartagena, C/ Doctor Fleming sn, 30202, cartagena, spain}
\email{jose.canovas@upct.es}
\address{ Houssem Eddine Rezgui$^{\text{2}}$, University of Carthage, faculty
of sciences of bizerte, (UR17ES21), "Dynamical Systems and their Applications"
\\ 7021, Jarzouna, Tunisia}
\email{rezguihoussemeddinne@gmail.com}

\date{}

\begin{abstract}
We consider the logistic family and apply the $q$-deformation $\phi_q(x)=\frac{1-q^x}{1-q}$. We study the stability regions of the fixed points of the $q$-deformed logistic map and the regions where the dynamic is complex through topological entropy and Lyapunov exponents. Our results show that the dynamic of this deformed family is richer than that of the $q$-deformed family studied in \cite{29}.
\end{abstract}


\maketitle

\section{Introduction}

It is known that $q$-deformations appeared in mechanic and thermodynamic systems in the so-called nonextensive statistical mechanics \cite{tsallis,4,3}, and they produced good results in quantum theory \cite{1,2}. However, this paper deals with its applications to discrete dynamical systems. It is remarkable that, recently, $q$-deformations and $q$-deformed maps have received the attention of many researchers \cite{9,7,29,12,13,3,5,4,10,8,6}. Let us introduce the general framework of this paper.

Recall that a discrete system is a pair $(X,f_{n})$, where $(X,d)$ is a (compact) metric space and $(f_{n})_{n\in \N}\subset \mathcal{F}$ where $\mathcal{F}=\{f:X\to X, f$ is continuous$\}$. For $x\in X$, its orbit is giving by the recursive sequence $x_{n+1}=f_{n}(x_{n})$, with initial condition $x_{0}=x$. If $f_{n}=f$ is constant, then the pair $(X,f)$ is a classical discrete dynamical system. If the sequence $f_{n}$ is periodic, we have a discrete periodic system. In this paper, we assume that $X=I=[0,1]$.

Several researchers have discussed the $q$-deformation of a function $f$, which implies the addition of a new parameter $q$ to the discrete dynamical system in such a way that $q$ tends to $1$, we recover the seminal discrete dynamical system generated by the function $f$. A $q$-deformation is a family of orientation preserving homeomorphisms of $I=[0,1]$. There are several families of $q$-deformations \cite{9,gupta}, which can be applied to different discrete dynamical systems (see e.g., \cite{9,12,8}).

In this paper, following \cite{9,gupta}, we consider the $q$-deformation given by
$$
\phi_{q}(x)=\frac{1-q^{x}}{1-q},
$$
where $x\in I$ and $q\in (0,+\infty)$, to produce a two-parameter family of $q$-deformed logistic maps. Recall that the logistic family is given by
$$
f_r(x)=rx(1-x),
$$
where $x\in I$ and $r\in (0,4]$. This two-parameter family is given by $Q_{r,q}=f_{r}\circ \phi_{q}$. Note that $Q_{r,q}$ is the two power of the periodic sequence
$$
(\phi_{q},f_{r},\phi_{q},f_{r},\dots),
$$
and applying the displacement operator to this sequence we get the shifted periodic sequence
$$
(f_{r},\phi_{q},f_{r},\phi_{q},\dots),
$$
from which we get the map $\Phi_{q,r}=\phi_{q}\circ f_{r}$. Both maps, $Q_{r,q}$ and $\Phi_{q,r}$ share most of their dynamical properties (see \cite{11}). Thus, it is more suitable working with $\Phi_{q,r}$ instead of $Q_{r,q}$.

In \cite{gupta}, the authors studied the dynamics of the map $Q_{r,q}$. They study the doubling period cascade to chaos and characterize the region where topological entropy is positive for $q\in (0,1)$. In this paper, following the ideas from \cite{29,12}, we analyze the dynamics of $\Phi_{q,r}$ by studying the stability of equilibrium points and characterizing the parameter region where its dynamic is chaotic. To do this characterization, we compute its topological entropy with fixed accuracy and estimate its Lyapunov exponents for a bigger parameter region, proving that interesting phenomena also appear for $q>1$.

In addition, $q$-deformed logistic maps exhibit the so-called Parrondo's paradox, where the composition of two dynamically simple maps produces a chaotic dynamical system \cite{11,12,14,gupta}. They can be considered an extremal example of this paradox since the $q$-deformations are homeomorphisms of $I$, which never can produce complicated dynamics. Therefore, we think these models are helpful for better understanding this paradox. Remarkably, this paradox has applications in physics, social sciences, and biology \cite{15,16,17,18}. We will show that the parameter region $q<1$ is highly suitable to produce the paradox.

On the other hand, it is then natural to ask whether about the dynamics of a deformed map when several $q$-deformations are  applied. We study the behavior of the periodic system $$(\phi_{q_{1}},\dots,\phi_{q_{k}},f_{r},\phi_{q_{1}},\dots, \phi_{q_{k}},f_{r},\dots),$$ for some $k\in\N$. The deformed system defined by the map $Q_{r,q_{k},\dots,q_{1}}=f_{r}\circ \phi_{q_{k},\dots,q_{1}}$. But, as in the case of a single $q$-deformation, for practical reasons we can make use of $\Phi_{q_{k},\dots,q_{1},r}= \phi_{q_{k},\dots,q_{1}}\circ f_{r}$. These generalized $q$-deformed logistic maps were studied in \cite{29} under a different $q$-deformation. A similar analysis for the $q$-deformation of this paper is not possible. We will point out the differences between the two generalized $q$-deformations.

This paper is organized as follows. In Section \ref{sec2}, we present some definitions and results which are useful for a proper understanding of the results presented in Sections \ref{sec4} and \ref{sec5}. Next, in Section \ref{sec3}, we study some basic properties of the $q$-deformed and $q$-deformations analyzed in this paper. Section \ref{sec4} is devoted to analyzing the stability of the equilibrium points of $\Phi_{q,r}$, while Section \ref{sec5} contains our result in characterizing the chaotic dynamics. The paper finishes with a section on conclusions and future work.

\section{ Definitions and preliminary results.\label{sec2}}

\subsection{Attractors and negative Schwarzian derivative.}

To state our results, we first introduce some notations (see, e.g. \cite{19}). Let $I=[0,1]$ and let $f:I \to I$ be a continuous interval map. Let $\mathbb{Z},\ \mathbb{Z}_{+}$, and $\mathbb{N}$ be the sets of integers, non-negative integers, and positive integers.

For $n\in \mathbb{Z_{+}}$ denote by  $f^{n}$ the $n$-$\textrm{th}$ iterate of $f$ i.e., $f^{0}$=identity, $f^1=f$, and $f^{n}=f\circ f^{n-1}$.
For any $x\in X$, the subset $O_{f}(x)= \{f^{n}(x): n\in\mathbb{Z}_{+}\}$ is called the \textit{orbit} of $x$ (under $f$). The accumulation points of an orbit, denote by $\omega(x,f)$, is the $\omega$-limit set of $x$ under $f$.
A set $A\subset X$ is called \textit{invariant} by $f$ (resp. \textit{strictly invariant by $f$}) if $f(A)\subset A$ (resp. $f(A)=A$). The basin of attraction $B(A)$ is giving by $$B(A)=\{x: \omega(x,f)\subset A\}.$$

A forward invariant compact set $M$ is called a (\textit{minimal}) \textit{metric attractor} if $B(M)$ has a positive Lebesgue measure, and if
$M^{'}$ is a strictly invariant compact set strictly contained in $M$, then $B(M^{'})$ has zero Lebesgue measure (see \cite{22, MI}).

A point $x_0\in I$ is called \textit{periodic} of period $n\in\mathbb{N}$ if $f^{n}(x_0)=x_0$ and $f^{i}(x_0)\neq x_0$ for $1\leq i\leq n-1$. In particular, if $n = 1$, then $x_0$ is called a \textit{fixed point} of $f$. According to \cite{EL}, a fixed point $x_0$ is said to be:
\begin{enumerate}
  \item Locally stable if for each $\varepsilon>0$ there exists a neighborhood $V$ of $x_0$ such that for any $x\in V$, then $|f^{n}(x)-x_0|<\varepsilon$, for each $n\in\N$.
  \item Attracting if there exists a neighborhood $V$ of $x_0$ such that $\displaystyle\lim_{n\to+\infty} f^{n}(x)=x_0$ for each $x\in V$.
  \item Locally asymptotically stable (LAS) if it is both, locally stable and attracting.
  \item Globally asymptotically stable (GAS) if it is LAS and $V=\mathrm{int}(I)$ where $\mathrm{int}(I)=(0,1)$ denotes the interior of $I$.
\end{enumerate}

Furthermore, if a fixed point $x_0$ is LAS (resp. GAS), then $|f^{'}(c)|\leq1$ whenever $f$ is differentiable.

If $f:I\to I$ is differentiable, a point $c$ is called a \textit{turning point} if $f^{'}(c)=0$. When $f$ is differentiable enough, a turning point $c$ is \textit{nonflat} if $f^{(n)}(c)\neq 0$ for some $n\geq 2$.
 If the map $f:I\to I$ is $\mathcal{C}^{3}$, we defined its \textit{Schwarzian derivative} at $x$ by
 $$
 \mathcal{S}(f)(x)=\frac{f^{'''}(x)}{f^{'}(x)}-\frac{3}{2}\left( \frac{f^{''}(x)}{f^{'}(x)} \right) ^{2}
 $$
 for all $x\in I$ such that $f^{'}(x)\neq 0$. The Schwarzian derivative of composition is obtained as follows (see \cite[Theorem 2.1]{20}).

\begin{thm}\cite{20}.\label{zz}
Let $f$ and $g$ be $\mathcal{C}^{3}$-maps, then $$\mathcal{S}(f\circ g)(x)=(\mathcal{S}(f)(g(x)))\cdot (g^{'}(x))^{2}+ \mathcal{S}(g)(x),$$ at those points $x\in I$ for which these expressions are well defined.
\end{thm}

Now we recall that a map $f:I\to I$ is unimodal if $f(0)=f(1)=0$ and if there exists a unique turning point $c\in (0,1)$ such that $f$ is strictly increasing in $[0,c)$ and strictly decreasing in $(c,1]$. Note that if $f$ is differentiable, then $f'(c)=0$. An interval $J\subset I$ is called wandering if the intervals $\{J,f(J),\dots,f^{n}(J),\dots\}$ are pairwise disjoint and $J$ is not contained in the basin of a LAS periodic orbit of $f$. As we will see below, the following result establishes the connection between the Schwarzian derivative and the existence of wandering intervals.

\begin{thm}\label{ww}\cite[Chapter II, Theorem 6.3]{21}.
Let $f:I\to I$ be a $\mathcal{C}^{3}$ unimodal map with negative Schwarzian derivative and such that $f^{''}(c)\neq0$ at the unique critical point $c$ of $f$. Then $f$ has no wandering intervals.
\end{thm}

On the other hand, it is interesting to know whether a LAS fixed point is also GAS, as the next result shows.

\begin{thm}\cite{20}.\label{liz}
Let $f:I\to I$ be a $\mathcal{C}^{3}$ unimodal map with negative Schwarzian derivative and it has a unique fixed point $c>0$ which is LAS, then it is GAS.
\end{thm}

For unimodal maps, the number of LAS periodic orbits are bounded as follows.

\begin{thm}\cite[Theorem 2.7]{20}.
Let $f:I\to I$ be a $\mathcal{C}^{3}$ unimodal map with negative Schwarzian derivative. Then each LAS periodic orbit attracts at least one critical point or belongs to the boundary of $I$.
\end{thm}

According to  \cite[Corollary 1]{22}, the structure of the metric attractors of a good enough unimodal map is described as follow

\begin{thm}\cite{22}.\label{aa}
Let $f:I\to I$ be a $\mathcal{C}^{3}$ unimodal map, then every metric attractor of $f$ is either
\begin{enumerate}
  \item a periodic orbit;
  \item a finite union of pairwise disjoint subintervals $I_{1},I_{2},\dots,I_{k}$ such that $f^{k}(I_{i})=I_{i}$ and $f^{k}_{|_{I_{i}}}$ has a dense orbit for $i\in\{1,\dots,k\}$;
  \item a Cantor set,
\end{enumerate}
and there is at most one metric attractor of a type different from type (1).
\end{thm}

Moreover, if the map $f$ has a negative Schwarzian derivative the number of attractors is at most two and just one if zero is not LAS (see e.g., cite{21}).

\subsection{Topological and physically observable chaos.}

Adler, Konheim and McAndrew introduced the topological entropy $h(f)$ of a continuous map $f$ on compact topological spaces \cite{23}. When $f:I\to I$ is unimodal, and in fact for piecewise monotone maps, Misiurewicz and Szlenk proved a simple version to define topological entropy as follows.

\begin{thm}\cite[Theorem 1]{24}
Let $f:I\to I$ be a piecewise monotone interval map, in particular unimodal map. Then, the topological entropy of $f$ is given by $$h(f)=\lim_{n\to+\infty}\frac{1}{n}\log c_{n},$$ where $c_{n}$ denotes the number of pieces of monotonicity of $f^{n}$.
\end{thm}

For unimodal maps, it is possible to compute the topological entropy practically. Although different algorithms have been proposed to calculate topological entropy, we will use the one proposed by Block et al. in \cite{28}. This algorithm compares the topological entropy of a unimodal map with that of the family of tent maps. Recall that tent maps are strictly piecewise linear maps with fixed slopes $s$ or $-s$, $1\leq s\leq 2$, whose topological entropy exactly equal to $\log (s)$. The algorithm, briefly described in \cite{12} as a slight modification of that from \cite{28}, is based on the classical bisection method to find the roots of real equations.

Let $f:I\to I$ be a continuous map. A pair of two different points $(x,y)\in I^{2}$ is \textit{proximal} if $$\displaystyle\liminf_{n\to+\infty}d(f^{n}(x),f^{n}(y))=0,$$ \textit{asymptotic} if $$\displaystyle\limsup_{n\to+\infty}d(f^{n}(x),f^{n}(y))=0,$$
and \textit{scrambled} or \textit{Li-Yorke} if it is proximal but not asymptotic.
A set $S\subseteq I$ is \textit{LY-scrambled} for $f$ if it contains at least two distinct points and every pair of distinct points in $S$ is scrambled. We say that $f:I\to I$ is LY-chaotic if there exists an uncountable LY-scrambled set \cite{25}. Moreover, it is known that a map with positive topological entropy is chaotic in the sense of Li and Yorke \cite{26}. There exist maps with topological entropy equal to zero and LY-chaotic (see \cite{27}), but they have to have wandering intervals \cite{baljim}. On the other hand, it is known that positive topological entropy interval maps may have just periodic attractors, and so chaos is not physically observable \cite{21,thunberg}.

For differentiable unimodal maps, the Lyapunov exponent of $f$ at $x\in I$ (see \cite{21}), is defined as
$$
\lambda(f,x)=\displaystyle \limsup_{n\to+\infty} \frac{1}{n} \sum_{i=1}^{n}\log|f^{'}(f^{i}(x))|.
$$
If we compute it at the turning point $c$, we can see that if the map $f$ has an attractor of type (1) according to Theorem \ref{aa}, then $\lambda(f,c)<0$. So, it is necessary to have a positive Lyapunov exponent to present observable chaos.

\section{Basic properties of $q$-deformations\label{sec3}}

First, the map $\Phi_{q,r}=\phi_{q}\circ f_{r}$ shares most of its dynamic properties with $Q_{r,q}$ (see \cite{11}). In particular, they have the same topological entropy. Besides, their Lyapunov exponents at the turning point have the same sign. Furthermore, it is easy to show that $\Phi_{q,r}$ verifies the following properties:

\begin{enumerate}
\item[(P1)] It follows from the chain rule and the fact that $\phi_{q}$ is a homeomorphism that the map  $\Phi_{q,r}$ is unimodal. It has a unique turning point $x_{0}=\frac{1}{2}$. To check this fact it suffices to realize that
$$
\Phi_{q,r}'(x)=\phi_q'(f_r(x))f_r'(x),
$$
and $\phi_q'(f_r(x))\neq 0$ for all $x\in [0,1]$.

\item[(P2)] We show that $x_0=\frac{1}{2}$ is non-flat because
\begin{eqnarray*}
\Phi_{q,r}''\left( \frac{1}{2} \right) &=& \phi_{q}''\left( f_{r}\left( \frac{1}{2}\right) \right) \cdot f_r'\left( \frac{1}{2}\right) ^2 +\phi_{q}'\left( f_{r}\left( \frac{1}{2}\right) \right)\cdot f_{r}^{''}\left( \frac{1}{2} \right) \\
&=& -\frac{2rq^{r/4}\log q}{1-q},
\end{eqnarray*}
if $q\neq 1$, and clearly, $q^{r/4}\log q\neq 0$ for all $q\neq 1$.

\item[(P3)] A straightforward computation shows that $\mathcal{S}(\phi_{q})=-\frac{(\log q)^2}{2}$. Thus, by Theorem \ref{zz}, the Schwarzian derivative
\begin{eqnarray*}
\mathcal{S}(\Phi_{q,r})&=&\mathcal{S}(\phi_q)(f_r(x))\cdot f'_r(x)^2 + \mathcal{S}(f_{r}) \\
&=& -(r(1-2x))^2\frac{\log ^2q}{2} -\frac{6}{(1-2x)^{2}}<0.
\end{eqnarray*}

\end{enumerate}

Note that when several $q$-deformations are applied, the properties (P1)-(P3) hold for the map $\Phi _{q_k,...,q_1,r}:=\phi _{q_k}\circ ...\circ \phi _{q_1}\circ f_{r}$. (P1) and (P3) follow inductively in a similar way. To check (P2) we use induction as well, assume that $\Phi_{q_k,...,q_1,r}''\left( \frac{1}{2} \right) \neq 0$, and compute
\begin{eqnarray*}
\Phi_{q_{k+1},q_k,...,q_1,r}'' \left( \frac{1}{2} \right) &=&
\phi_{q_{k+1}}'' \left( \Phi_{q_k,...,q_1,r} \left( \frac{1}{2} \right)  \right)\cdot
\left( \Phi_{q_k,...,q_1,r}'\left( \frac{1}{2} \right) \right) ^2   \\
&+& \phi_{q_{k+1}}' \left( \Phi_{q_k,...,q_1,r} \left( \frac{1}{2} \right)  \right) \cdot \Phi_{q_k,...,q_1,r}''\left( \frac{1}{2} \right) \\
&=& -\frac{q_{k+1}^{\Phi_{q_k,...,q_1,r} \left( \frac{1}{2} \right)}\log q_{k+1}}{1-q_{k+1}} \cdot \Phi_{q_k,...,q_1,r}''\left( \frac{1}{2} \right) \neq 0,
\end{eqnarray*}
where we have applied that $\Phi_{q_k,...,q_1,r}'\left( \frac{1}{2} \right)=0$.

Regarding the map $\phi_q$ we have the following property.

\begin{Prop}
For $0<q_1<q_2<1<q_3<q_4$ we have that $\phi _{q_1}(x)>\phi _{q_2}(x)>x>\phi _{q_3}(x)>\phi _{q_4}(x)$ for all $x\in (0,1)$.
\end{Prop}

\textbf{Proof. }
Fix $q_1,q_2\in (0,\infty)\setminus \{1\}$. Define $F_{q_1,q_2}(x):=\phi_{q_1}(x)-\phi_{q_2}(x)$ for all $x\in [0,1]$. Note that $F_{q_1,q_2}(0)=F_{q_1,q_2}(1)=0$. We claim that there is not $x^*\in (0,1)$ such that $F_{q_1,q_2}(x^*)=0$. Otherwise, there must exist two points $x_1,x_2\in (0,1)$ such that $F_{q_1,q_2}'(x_i)=0$ for $i=1,2$. However, a straightforward computation shows that
$$
\frac{\partial F_{q_1,q_2} (x)}{\partial x}=\frac{q_2^x\log q_2}{1-q2}-\frac{q_1^x\log q_1}{1-q1},
$$
and thus, the equation $F_{q_1,q_2}'(x)=0$ has a unique solution
$$
x=-\frac{\log \left( \frac{(1 - q_2) \log q_1}{(1 - q_1) \log q_2} \right) } {\log q_1 - \log q_2},
$$
and therefore $x^*$ cannot exist.
We obtain a similar result in a similar way for $q_1\in (0,\infty)\setminus \{1\}$ and $q_2=1$ by taking the map $F_{q_1}(x):=\phi_{q_1}(x)-x$.

Next, we fix $x=\frac{1}{2}$ and consider the map
$$
\varphi (q):=\phi_q(1/2)=\frac{1-\sqrt{q}}{1-q}
$$
if $q\neq 1$, and $\varphi (1)=1/2$. Note that this map can be rewritten as
$$
\varphi (q)=\frac{1}{1+\sqrt{q}},
$$
which is strictly decreasing. This finishes the proof.
$\square$

\begin{figure}[htbp]
\begin{center}
\begin{tabular}{cc}
(a) \includegraphics[width=0.35\textwidth]{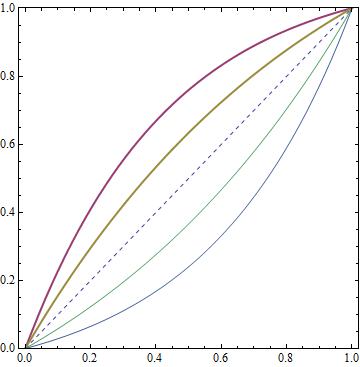} &
(b) \includegraphics[width=0.35\textwidth]{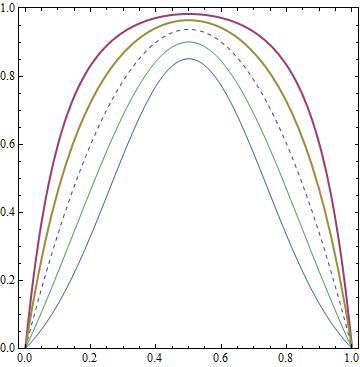} \\

\end{tabular}
\caption{(a) Graph of the identity (dashed line), $\phi _{1/3}$ and $\phi _{1/10}$ (thick lines) and  $\phi _3$ and $\phi _{10}$ (thin lines).(b) Graph of the maps $f_{3.75}$ (dashed line), $\Phi _{1/3,3.75}$ and $\Phi _{1/10,3.75}$ (thick lines) and $\Phi _{3,3.75}$ and $\Phi _{10,3,75}$ (thin lines). }
\label{graficas}
\end{center}
\end{figure}

Figure \ref{graficas} shows the graph of the maps $\phi_q$ and $\Phi_{q,3.75}$ for several values of $q$. As a consequence, we obtain the following result whose proof is immediate.

\begin{Prop}
For $0<q_1<q_2<1<q_3<q_4$ and for $r\in (0,4]$, we have that $\Phi _{q_1,r}(x)>\Phi _{q_2,r}(x)>f_r(x)>\Phi _{q_3,r}(x)>\Phi _{q_4,r}(x)$ for all $x\in (0,1)$.
\end{Prop}

\section{Stability of fixed points\label{sec4}}

\subsection{One $q$-deformation applied}
In this section we deal with the stability of the fixed points of $\Phi_{q,r}$. First, note that
$$
\Phi_{q,r}(0)=\phi_q(0)=f_r(0)=0,
$$
and thus $0$ is a fixed point for $\Phi_{q,r}$, $\phi_q$ and $f_r$. Finding non-zero fixed points requires to solve the equation
$$
x=\Phi_{q,r}(x)=\frac{1-q^{rx(1-x)}}{1-q},
$$
which cannot be solved in a closed way when $q\neq 1$ and only numerical methods can be used to obtain an approximate solution. Recall that $f_{r}$ ($q=1$) has a non-zero fixed point $x^{*}=\frac{r-1}{r}$ whenever $r>1$. It is well-known that $0$ is asymptotically stable if $0<r\leq 1$ and $x^{*}$ is asymptotically stable if $1<r\leq3$.

Regarding the fixed point $0$ we can prove the following result.

\begin{thm}\label{el0}
The fixed point $0$ is LAS when $(q,r)\in \mathcal{S}_0$, where
\begin{eqnarray*}
\mathcal{S}_0 &:=& \left\{ (q,r)\in (0,\infty)\times (0,4]: \frac{r \log(q)}{q-1}<1 \right\} \\
&& \cup \left\{ (q,r)\in (0,\infty)\times (0,4]: r \log(q)=q-1 \mbox{ and }q\leq 3  \right\} .
\end{eqnarray*}
\end{thm}

\textbf{Proof.}
Fix $q\neq 1$ and note that
$$
\Phi_{q,r}'(0)=(f_{r}'(0)\cdot \phi_{q})'(0)=\frac{r \log(q)}{q-1}>0.
$$
Clearly, if $\frac{r \log(q)}{q-1}<1$, the fixed point $0$ is LAS. If $\frac{r \log(q)}{q-1}=1$, by \cite[Theorem 1.5]{EL}, the fixed point $0$ is LAS in the following cases:
\begin{itemize}
\item When $\Phi_{q,r}''(0)<0$, which reads as
$$
r\log q=q-1<2,
$$
which is equivalent to $q<3$.

\item When $q= 3$, $\Phi^{'''}_{q,r}(0)<0$, which is fulfilled because $\Phi^{'''}_{q,r}(0)=-8<0$.
\end{itemize}
which concludes the proof.
$\square$

\begin{figure}[htbp]
\begin{center}

 \includegraphics[width=0.45\textwidth]{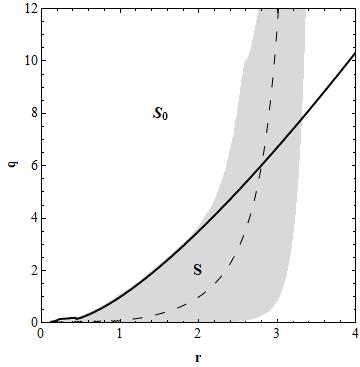}
\caption{Graph of the curve $\Phi_{q,r}(1/2)=1/2$ (dashed line), graph of the curve $r \log(q)=q-1$ (thick line), $S_0$ the region where $0$ is LAS and $S$ the region where a nonzero fixed point is LAS.}
\label{fixedpoints1}
\end{center}
\end{figure}

Figure \ref{fixedpoints1} shows the shape of the set $\mathcal{S}_0$. Next, we analyze the existence of non-zero fixed points. We distinguish two cases:

\begin{enumerate}
\item There is $\varepsilon >0$ such that $\Phi_{q,r}(x)>x$ for all $x\in (0,\varepsilon )$. As $\Phi_{q,r}(1)=0$, then there exists another fixed point. Note that $0$ is not LAS. As the Schwarzian derivative is negative, there is just one metric attractor.

\item There is $\varepsilon >0$ such that $\Phi_{q,r}(x)<x$ for all $x\in (0,\varepsilon )$. Note that $0$ is not LAS. As the Schwarzian derivative is negative, two metric attractors can coexist.
    Here we can have 0, 1, or 2 non-zero fixed points, and only one of them can be LAS.

\end{enumerate}

Then, we can obtain the following result.

\begin{thm}\label{lasgas}
Let $\Phi_{q,r}$ be the $q$-deformed logistic map. Then the following statements hold.
\begin{enumerate}
\item If $0$ is LAS, then it is GAS if and only if $f(x)<x$.
\item If there exists a non-zero LAS fixed point, then it is GAS if and only if $0$ is not LAS.
\end{enumerate}
\end{thm}

\textbf{Proof.}
It is straightforward taking into account that $\Phi_{q,r}$ has negative Schwarzian derivative and Theorem \ref{liz}.
$\square$

It is not possible to obtain the non-zero fixed points of $\Phi_{q,r}$ in an analytical way. So, we must obtain them numerically and get the approximate value of the first derivative at such points. The line $\Phi_{q,r}(1/2)=1/2$ provides some control on our computations since the inequality $\Phi_{q,r}(1/2)>1/2$ ensures the existence of at least one non-zero fixed point. Figure \ref{fixedpoints1} shows the shape of the region $\mathcal{S}$ where a non-zero fixed point is LAS from our numerical simulations.

\begin{figure}[htbp]
\begin{center}
\begin{tabular}{cc}
(a) \includegraphics[width=0.35\textwidth]{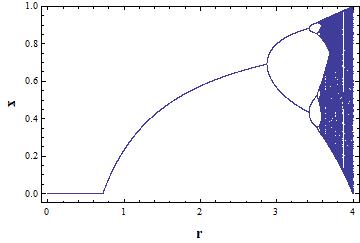} &
(b) \includegraphics[width=0.35\textwidth]{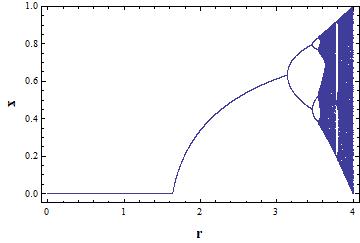} \\

\end{tabular}
\caption{(a) Bifurcation diagram of $\Phi_{0.5,r}$ for $r\in (0,4]$. (a) Bifurcation diagram of $\Phi_{2.5,r}$ for $r\in (0,4]$. We draw the last 100 iterates of orbits of length 10000 and initial condition $x=0.45$.}
\label{bif1q}
\end{center}
\end{figure}

Note that the thick line of Figure \ref{fixedpoints1} is given by $\Phi_{q,r}'(0)=1$, and it seems to be the frontier of the region $\mathcal{S}$ for $q<3$ (compare with Theorem \ref{el0}). This fact indicates that there is a saddle-node bifurcation (see \cite{kutne}). Figure \ref{bif1q} shows several bifurcation diagrams that show this fact. The bifurcation diagrams also show the existence of period-doubling bifurcation when the non-zero fixed point destabilizes with derivative $-1$. The well-known period-doubling route to chaos can also be appreciated.

However, according to \cite{kutne}, the map $\Phi_{q,r}$ does not hold the necessary conditions to guarantee that $0$ destabilizes following a saddle-node bifurcation, because although
$$
\frac{\partial ^2}{\partial x^2}\Phi_{q,r}(0)=q-3\neq 0
$$
for $q<3$, we have by straightforward calculation that both
\begin{equation}
\frac{\partial }{\partial r}\Phi_{q,r}(0)=\frac{\partial }{\partial q}\Phi_{q,r}(0)=0. \label{nondeg}
\end{equation}
It is worth mentioning that the condition (\ref{nondeg}) makes possible richer bifurcation scenarios. For example, if both derivatives differed from zero, we would have a saddle-node bifurcation from a LAS fixed point $0$ to another LAS fixed point different from $0$. Numerical simulations show that this saddle-node bifurcation happens when $q<3$, but then the stability region of $0$ intersects the one of a non-zero LAS fixed point, making possible different bifurcations shown in this paper, as it is explained below.

\begin{figure}[htbp]
\begin{center}
\begin{tabular}{cc}
(a) \includegraphics[width=0.35\textwidth]{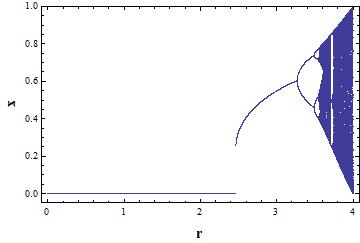} &
(b) \includegraphics[width=0.35\textwidth]{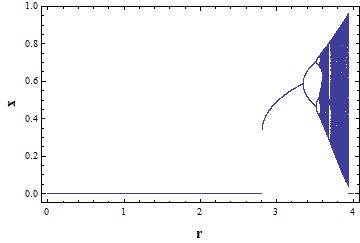} \\
(c) \includegraphics[width=0.35\textwidth]{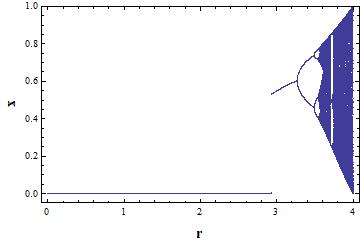} &
(d) \includegraphics[width=0.35\textwidth]{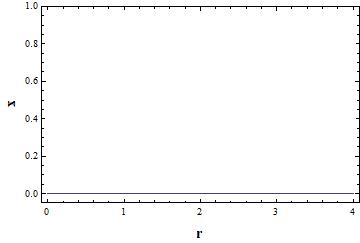} \\
\end{tabular}
\caption{(a) Bifurcation diagram of $\Phi_{6.5,r}$ for $r\in (0,4]$. (b) Bifurcation diagram of $\Phi_{12.5,r}$ for $r\in (0,4]$. We depict the last 100 iterates of orbits of length 10000 and initial condition $x=0.45$. (c) and (d) show the bifurcation diagrams when the initial condition is $x=0.001$.}
\label{bif2q}
\end{center}
\end{figure}

We denote by $\mathcal{S}_0(r):=\{ q\in (0,\infty):(q,r)\in \mathcal{S}_0\}$. Note that the map
$$
r_1(q):=\frac{q-1}{\log q}
$$
is strictly increasing with $\lim _{q\rightarrow 0}r_1(q)=0$ and  $\lim _{q\rightarrow \infty}r_1(q)=\infty$. So, $\mathcal{S}_0(r_1)\subset \mathcal{S}_0(r_2)$ for $r_1<r_2$. In particular, for $q>q_0=10.34665192905221$ we can check that $\mathcal{S}_0(r_1)=(0,4]$. As we will see, this will imply that the dynamics can be topologically chaotic but almost all the orbits are attracted by the fixed point $0$. On the other hand, the curve $\Phi_{q,r}(1/2)=1/2$ reads as
$$
r=r_2(q):=4\frac{\log \left( \frac{1+q}{2} \right) }{\log q},
$$
which is increasing and holds $\lim _{q\rightarrow 0}r_2(q)=0$ and  $\lim _{q\rightarrow \infty}r_1(q)=4$. In particular, we have that $r_1(q)=r_2(q)$ when $q=q_1=6.025724834504679$. Then, other possible bifurcations can be exhibit by $\Phi_{q,r}$ when the fixed point $0$ loses its stability. Figure \ref{bif2q} shows this fact.

This kind of bifurcation happens when two attractors coexist, and the orbit changes the basin of attraction. Here, we have that the fixed point $0$ is LAS and two additional fixed points $x_1<x_2$ also appear. The fixed point $x_1$ is not LAS, but it describes the basin of attraction of $0$, which contains $[0,x_1)\cup (x_1^*,1]$, where $x_1^*:=\max f^{-1}(x_1)$. Now, we have two possibilities: if the set $[x_1,x_1^*]$ is invariant by $f$, we will have a unique attractor on it. If it is not, that is, $\Phi_{q,r}(1/2)>x_1^*$, then almost all orbit inside $[x_1,x_1^*]$ escapes from it to the basin of attraction of the fixed point $0$. Figure \ref{bif2q} shows the bifurcation diagrams with two different initial conditions, one close enough to $0$ and the other close to the turning point.

\subsection{Several $q$-deformations applied\label{several}}

Here we consider the case when several $q$-deformations $\phi_{q_1},...,\phi_{q_k}$ are applied to $f_r$ to construct a new map $\Phi _{q_k,...,q_1,r}:= \phi_{q_k}\circ ...\circ \phi_{q_1}\circ f_r$. We find the main difference with the paper \cite{29}, because in that paper the composition of several $q$-deformation commuted and, even more, it was also a $q$-deformation. As a consequence, the study of the application of several $q$-deformations to $f_r$ can be reduced to the single case. In this paper, it is immediate to realized that neither the composition of $q$-deformations commutes, nor it is a new $q$-deformation. Then, studying the dynamics of the map $\Phi _{q_k,...,q_1,r}$ becomes harder. As in \cite{29}, we consider two different cases. First, we consider $q_1=...=q_k=q$ so that the map $\Phi _{q_k,...,q_1,r}=\Phi _{k;q,r}=\phi_q^k\circ f_r$. Next, we will analyze the case when  $q_1=...=q_k$ are chosen between two different values on $(0,\infty)$.

\subsubsection{The case $q_1=...=q_k=q$.} 
\begin{figure}[htbp]
\begin{center}
\begin{tabular}{ccc}
(a) \includegraphics[width=0.3\textwidth]{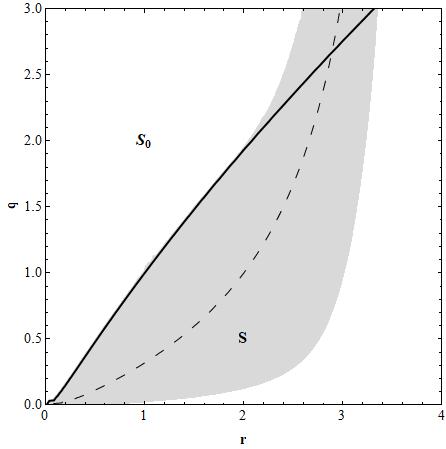} &
(b) \includegraphics[width=0.3\textwidth]{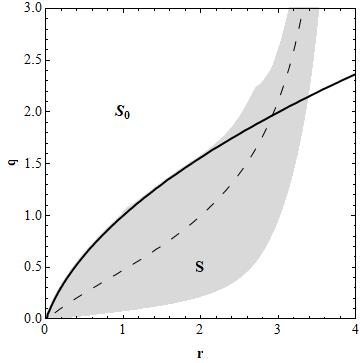} &
(c) \includegraphics[width=0.3\textwidth]{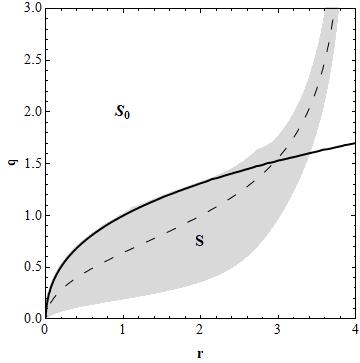}
\end{tabular}
\caption{For $k=2$(a), $k=3$(b) and $k=5$(c), graph of the curve $\Phi_{k;q,r}(1/2)=1/2$ (dashed line), graph of the curve $\Phi_{k;q,r}'(0)=1$ (thick line), $S_0$ the region where $0$ is LAS and $S$ the region where a nonzero fixed point is LAS. }
\label{stabk}
\end{center}
\end{figure}

We fix $k=2,3,5$ and check how the stability regions of fixed points evolve with $k$. As in the single case, we consider the curves $\Phi _{k;q,r}'(0)=1$ and $\Phi _{k;q,r}(1/2)=1/2$. Figure \ref{stabk} shows the stability region $\mathcal{S}_0$ of the fixed point $0$ and the numerically obtained stability region $\mathcal{S}$ of a non-zero fixed point. The pictures show that while $\mathcal{S}_0$ grows when $k$ grows, $\mathcal{S}$ seems to decrease.

\begin{figure}[htbp]
\begin{center}
\begin{tabular}{cc}
(a) \includegraphics[width=0.35\textwidth]{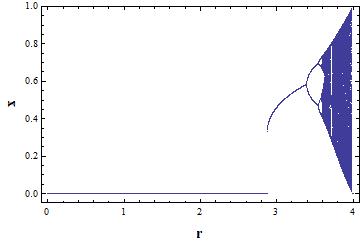} &
(b) \includegraphics[width=0.35\textwidth]{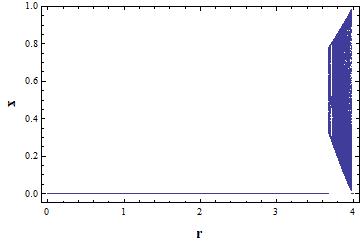} \\

\end{tabular}
\caption{Bifurcation diagrams of $\Phi_{2;3.6,r}$ for $r\in (0,4]$. We depict the last 100 iterates of orbits of length 10000 and initial condition $x=0.45$ (a) and $x=0.001$ (b).}
\label{bifk}
\end{center}
\end{figure}

Figure \ref{stabk} also shows that the region where saddle-node bifurcation happens is reduced, and the change of attractor type bifurcation appears for smaller values of $q$. For instance, in Figure \ref{bifk} we depict several bifurcation diagrams for $k=2$ showing this fact. It is straightforward to check that a similar condition to that of (\ref{nondeg}) happens in these cases, making possible richer bifurcation scenarios.

\subsubsection{Two different $q$-deformations applied.}
\begin{figure}[htbp]
\begin{center}
\begin{tabular}{cc}
(a) \includegraphics[width=0.35\textwidth]{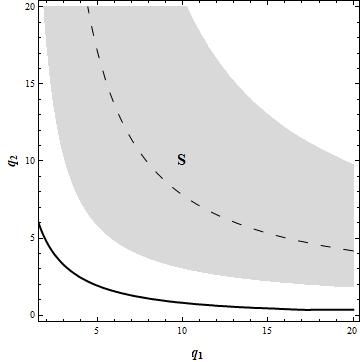} &
(b) \includegraphics[width=0.35\textwidth]{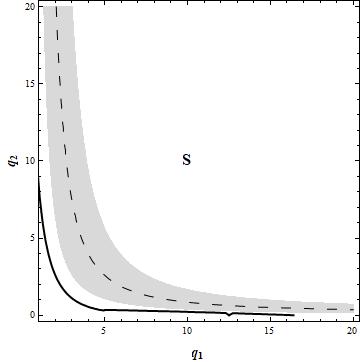} \\
(c) \includegraphics[width=0.35\textwidth]{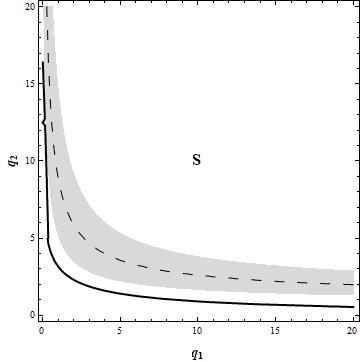} &
(d) \includegraphics[width=0.35\textwidth]{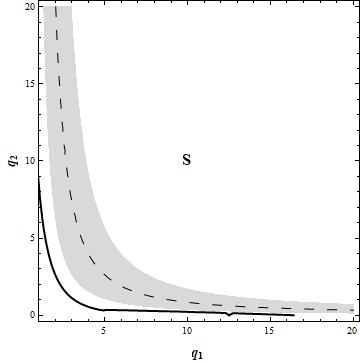} \\
\end{tabular}
\caption{Graph of the curve $\Phi(1/2)=1/2$ (dashed line), graph of the curve $\Phi '(0)=1$ (thick line), $S_0$ the region where $0$ is LAS and $S$ the region where a nonzero fixed point is LAS, for $\Phi$ equal to $\Phi_{q_1,q_2,3.5}$ (a), $\Phi_{q_1,q_1,q_2,3.5}$ (b), $\Phi_{q_1,q_2,q_2,3.5}$ (c) and $\Phi_{q_1,q_2,q_1,3.5}$ (d). }
\label{stab}
\end{center}
\end{figure}

Next, we consider the case where the map $\Phi _{q_k,...,q_1,r}$ is such that $q_k,...,q_1$ takes only two possible values at the same time. Then, the map depends on three parameters, which makes it quite difficult to present the results and graphics suitably. Hence, we are forced to fix some values of the parameter $r$ and check how the stability regions changes when more $q$-deformations are considered. In particular, we study the cases $\Phi_{q_1,q_2,r}$, $\Phi_{q_1,q_1,q_2,r}$, $\Phi_{q_1,q_2,q_2,r}$ and $\Phi_{q_1,q_2,q_1,r}$. Our results will show that the order in which the $q$-deformations are iterated influences the stability regions of the fixed point $0$ and the non-zero fixed point obtained numerically. This fact represents a huge difference with respect to the results of \cite{29}. To show it, we fix $r=3.5$ and depict the stability regions $\mathcal{S}_0$, where $0$ is LAS, and $\mathcal{S}$, where a non-zero fixed point is also LAS for the four maps mentioned above. The results are shown in Figure \ref{stab}.

\begin{figure}[htbp]
\begin{center}
\begin{tabular}{cc}
(a) \includegraphics[width=0.35\textwidth]{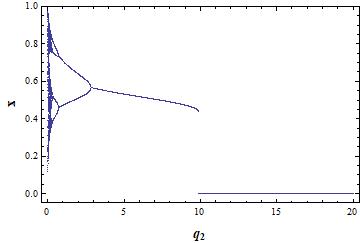} &
(b) \includegraphics[width=0.35\textwidth]{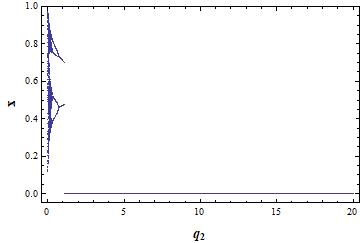} \\

\end{tabular}
\caption{Bifurcation diagrams of $\Phi_{3,q_2,q_2,3.5}$ for $q_2\in (0,20]$. We depict the last 100 iterates of orbits of length 10000 and initial condition $x=0.45$ (a) and $x=0.001$ (b).}
\label{bifstab}
\end{center}
\end{figure}

Figure \ref{stab} also shows that the destabilization of the fixed point $0$ need not happen via a saddle-node bifurcation but for a change of a basin of two coexistent attractors. Figure \ref{bifstab} shows the bifurcation diagram for $\Phi_{q_1,q_2,q_2,3.5}$.

\section{Topological and observable chaos\label{sec5}}

\subsection{One $q$-deformation applied}

\begin{figure}[htbp]
\begin{center}
\begin{tabular}{cc}
(a) \includegraphics[width=0.35\textwidth]{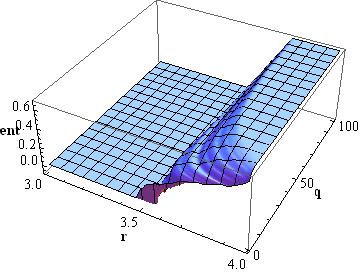} &
(b) \includegraphics[width=0.35\textwidth]{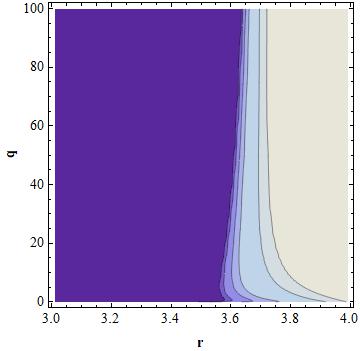} \\
(c) \includegraphics[width=0.35\textwidth]{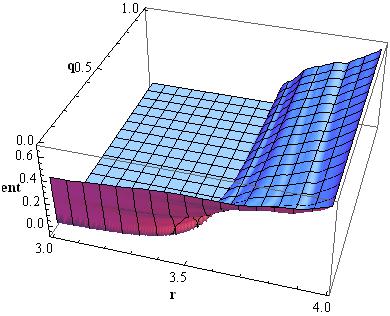} &
(d) \includegraphics[width=0.35\textwidth]{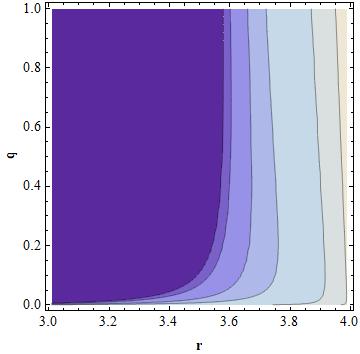} \\
(e) \includegraphics[width=0.35\textwidth]{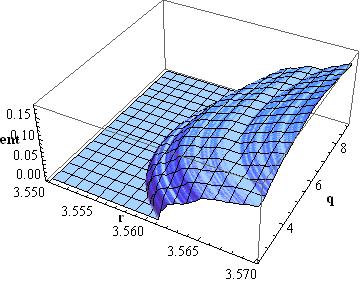} &
(f) \includegraphics[width=0.35\textwidth]{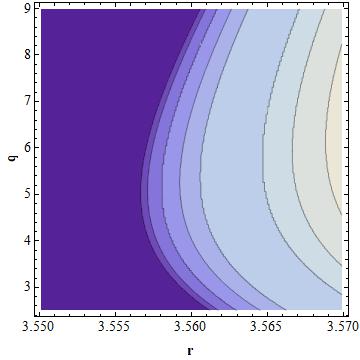}

\end{tabular}
\caption{Topological entropy of the $q$-deformed logistic map $\phi_q\circ f_r$ for $r\in [3,4]$ with step size $10^{-2}$ and $q\in [0,100)$ with step size $10^{-1}$ (a) and associated level curves (b), for $r\in [3,4]$ with step size $10^{-2}$ and $q\in [0,1)$ with step size $10^{-3}$ (c) and associated level curves (d), and for $r\in [3.55,3.57]$ with step size $10^{-3}$ and $q\in [2.5,9)$ with step size $10^{-2}$ (e) and associated level curves (f). The prescribed accuracy is $10^{-4}$. }
\label{fig:1}
\end{center}
\end{figure}

\begin{figure}[htbp]
\begin{center}
\begin{tabular}{cc}
(a) \includegraphics[width=0.35\textwidth]{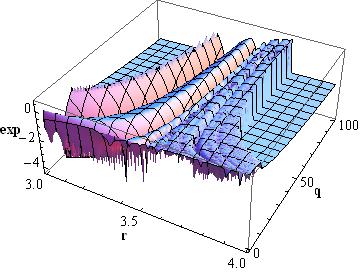} &
(b) \includegraphics[width=0.35\textwidth]{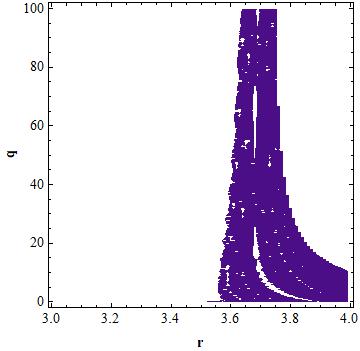}\\
(c) \includegraphics[width=0.35\textwidth]{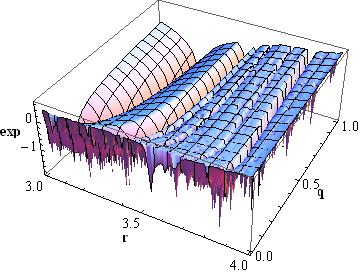} &
(d) \includegraphics[width=0.35\textwidth]{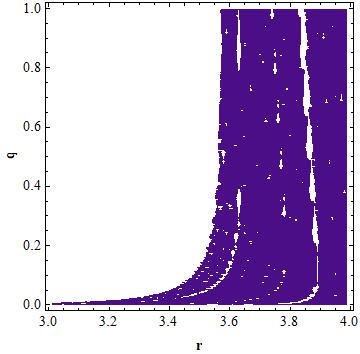}\\
(e) \includegraphics[width=0.35\textwidth]{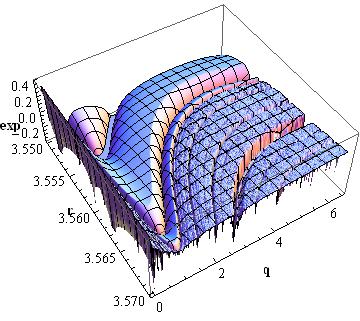} &
(f) \includegraphics[width=0.35\textwidth]{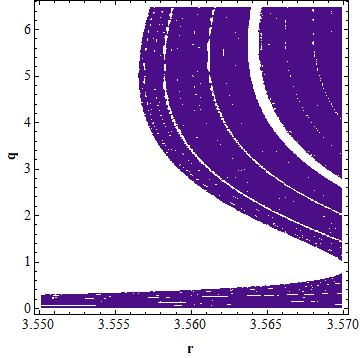}

\end{tabular}
\caption{(a) Estimations of the Lyapunov exponent at the turning point the $q$-deformed logistic map $\phi_q\circ f_r$ for $r\in [3,4]$ with step size $10^{-2}$ and $q\in [0,100)$ with step size $10^{-1}$ (b) Region where the estimations of the Lyapunov exponent is positive. (c)-(d) the same computations for $r\in [3,4]$ with step size $10^{-2}$ and $q\in [0,1)$ with step size $10^{-3}$.  (e)-(f) the same computations  for $r\in [3.55,3.57]$ with step size $10^{-3}$ and $q\in [2.5,9)$ with step size $10^{-2}$. We make the estimations using a orbit length of 10.000 points. }
\label{fig:5}
\end{center}
\end{figure}

By using the algorithm in \cite{28} performed as it is described in \cite{12}, we compute the topological entropy with an accuracy $10^{-4}$. In addition, we estimate the Lyapunov exponents to show the existence of observable chaos. The fact that all the maps have the same turning point simplifies the computation of topological entropy and Lyapunov exponents. Figures \ref{fig:1} and \ref{fig:5} show our computations.

We can see that the application of a $q$-deformation decreases the complexity of $\Phi_{q,r}$. However, the complexity remains for all $q>0$ since for $r=4$, we have that $\Phi_{q,4}([0,1/2])=\Phi_{q,4}([1/2,1])= [0,1]$ and hence we have a $2$-horseshoe, and the topological entropy is $\log 2$ (see e.g. \cite[Chapter 4]{30}). However, as Figure \ref{bif2q} points out, it may happen that, in this case, almost all the orbits converge to the fixed point $0$.

On the other hand, it is known that when $f_{r}$ has zero topological entropy when  $r\leq 3.5699....$. The $q$-deformation $\phi _q$ is always simple. Figure \ref{fig:1} shows that it is possible to have $h(\Phi_{q,r})>0$ for parameter values $r<3.5699...$. In fact, this is easily achieved for parameter values of $q$ smaller than 1. Yet, this is another example of the Parrondo's paradox "simple + simple= chaos". This paradox is exhibited for other $q$-deformations of the logistic map \cite{29,12}.

\subsection{Several $q$-deformations applied}
Now, we study a finite number of $q$-deformations applied to $f_{r}$ and analyze how the situation changes. Following the notation of Section \ref{several}, we consider two cases depending on the number of different $q$-deformations applied. As in the single case, we compute the topological entropy with accuracy $10^{-4}$ to prove the existence of topological chaos. To show that chaos can be physically observable, we estimate Lyapunov exponents.

\subsubsection{The case $q_1=...=q_k=q$.}
\begin{figure}[htbp]
\begin{center}
\begin{tabular}{cc}
(a) \includegraphics[width=0.35\textwidth]{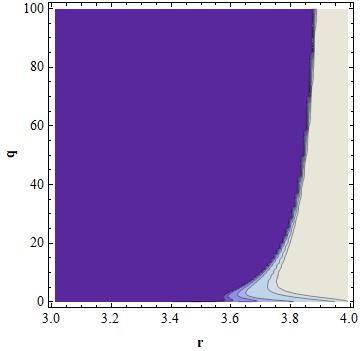} &
(b) \includegraphics[width=0.35\textwidth]{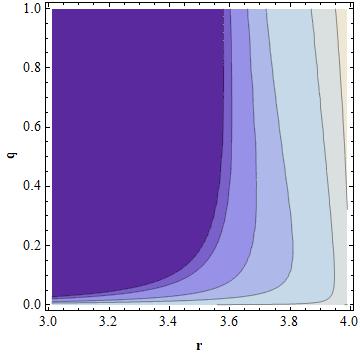} \\
(c) \includegraphics[width=0.35\textwidth]{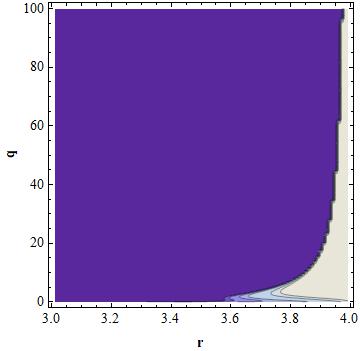} &
(d) \includegraphics[width=0.35\textwidth]{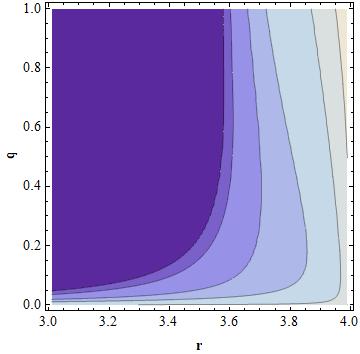} \\
(e) \includegraphics[width=0.35\textwidth]{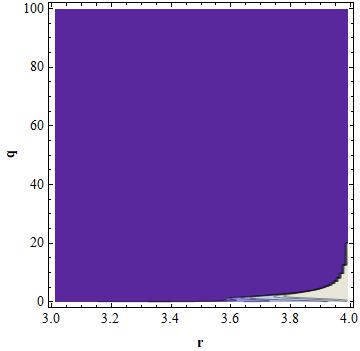} &
(f) \includegraphics[width=0.35\textwidth]{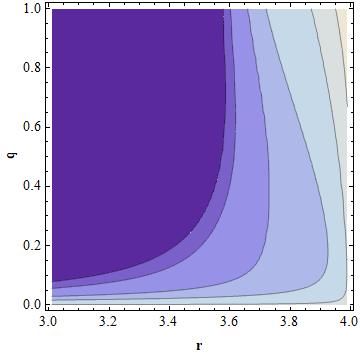}

\end{tabular}
\caption{Level curves of the topological entropy of $\Phi$ when $r\in [3,4]$ with step size $10^{-2}$ and $q\in [0,100)$ with step size $10^{-1}$ (left column) and $r\in [3,4]$ with step size $10^{-2}$ and $q\in [0,1)$ with step size $10^{-3}$ (right column) when $\Phi$ is equal to $\Phi_{2;q,r}$(a)-(b), $\Phi_{3;q,r}$(c)-(d) and $\Phi_{5;q,r}$(e)-(f). The prescribed accuracy is $10^{-4}$. }
\label{fig:11}
\end{center}
\end{figure}

\begin{figure}[htbp]
\begin{center}
\begin{tabular}{cc}
(a) \includegraphics[width=0.35\textwidth]{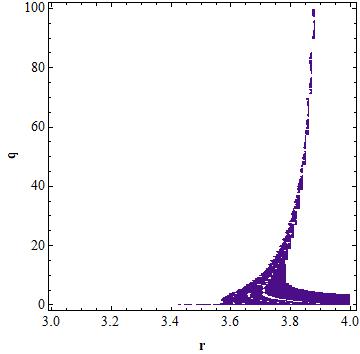} &
(b) \includegraphics[width=0.35\textwidth]{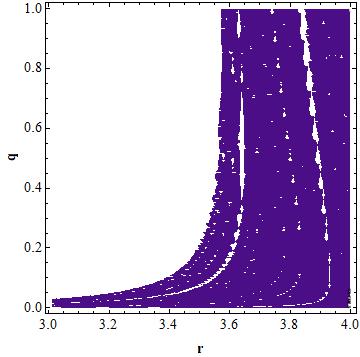} \\
(c) \includegraphics[width=0.35\textwidth]{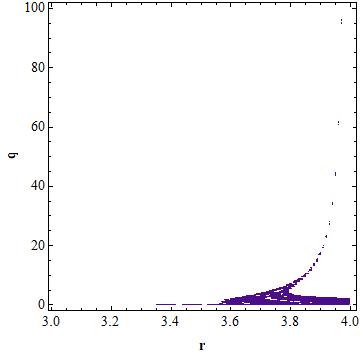} &
(d) \includegraphics[width=0.35\textwidth]{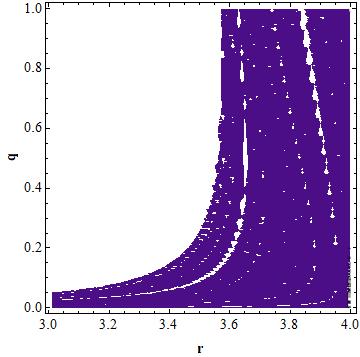} \\
(e) \includegraphics[width=0.35\textwidth]{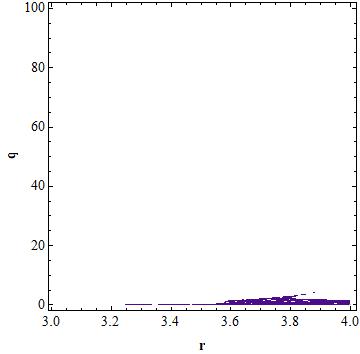} &
(f) \includegraphics[width=0.35\textwidth]{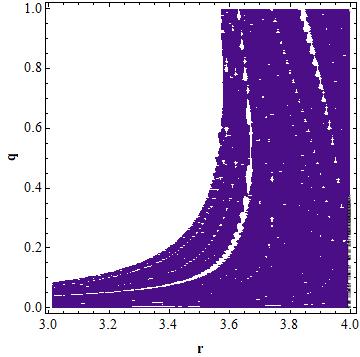}

\end{tabular}
\caption{Regions where the estimation of the Lyapunov exponent at the turning point of $\Phi$ is positive when $r\in [3,4]$ with step size $10^{-2}$ and $q\in [0,100)$ with step size $10^{-1}$ (left column) and $r\in [3,4]$ with step size $10^{-2}$ and $q\in [0,1)$ with step size $10^{-3}$ (right column) when $\Phi$ is equal to $\Phi_{2;q,r}$(a)-(b), $\Phi_{3;q,r}$(c)-(d) and $\Phi_{5;q,r}$(e)-(f). We make the estimations using a orbit of 10.000 points. }
\label{fig:14}
\end{center}
\end{figure}

Following the notation of Section \ref{several}, we compute the topological entropy of $\Phi_{k;q,r}$ for $k=2$, $3$ and $5$. The results are shown in Figure \ref{fig:11}. Similarly, we show the regions where the estimations of Lyapunov exponents are positive in Figure \ref{fig:14}. Note that we present our results on two different regions, the bigger one for $(q,r)\in (0,100]\times [3,4]$ presents a general overview of the situation and shows that the complexity decreases when both $q$ and $k$ increase. Next, we make a zoom in the region $(q,r)\in (0,1]\times [3,4]$ with different step sizes and show that the combination of simple dynamics may result in a complicated one. In this enlarged region, when $k$ increases, the region of parameter values exhibiting a chaotic behaviour may increase.

\subsubsection{Two different $q$-deformations applied.}

\begin{figure}[htbp]
\begin{center}
\begin{tabular}{cc}
(a) \includegraphics[width=0.35\textwidth]{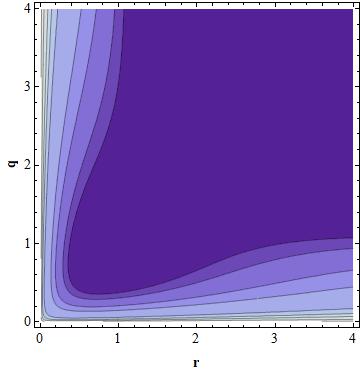} &
(b) \includegraphics[width=0.35\textwidth]{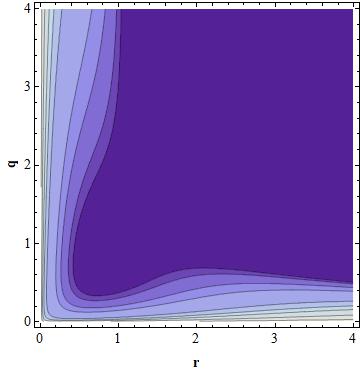} \\
(c) \includegraphics[width=0.35\textwidth]{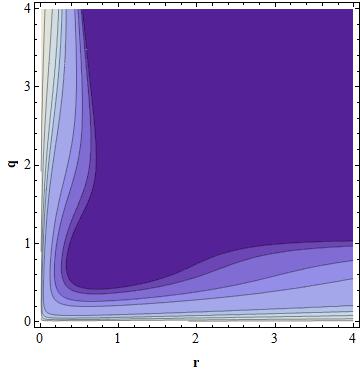} &
(d) \includegraphics[width=0.35\textwidth]{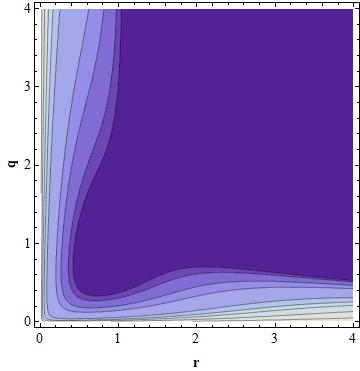} \\
\end{tabular}
\caption{Level curves of the topological entropy of the map $\Phi$ for $r=3.56$ and $q_1,q_2\in (0,4)$ with step size $10^{-1}$ for the maps $\Phi_{q_1,q_2,r}$(a), $\Phi_{q_1,q_1,q_2,r}$(b), $\Phi_{q_1,q_2,q_2,r}$(c) and $\Phi_{q_1,q_2,q_1,r}$(d). The accuracy of the computations is $10^{-4}$.  }
\label{top1}
\end{center}
\end{figure}

\begin{figure}[htbp]
\begin{center}
\begin{tabular}{cc}
(a) \includegraphics[width=0.35\textwidth]{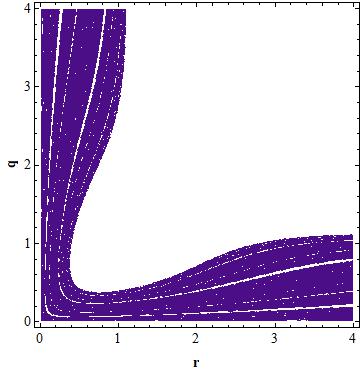} &
(b) \includegraphics[width=0.35\textwidth]{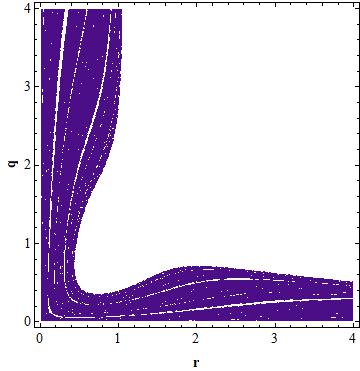} \\
(c) \includegraphics[width=0.35\textwidth]{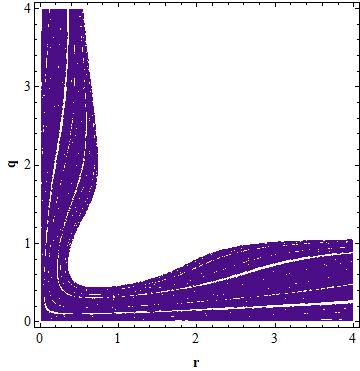} &
(d) \includegraphics[width=0.35\textwidth]{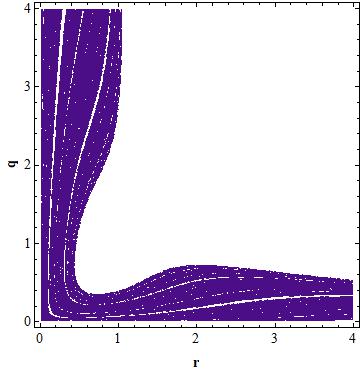} \\
\end{tabular}
\caption{Regions where the estimation of the Lyapunov exponent at the turning point of $\Phi$ is positive when $r=3.56$ and $q_1,q_2\in (0,4)$ with step size $10^{-1}$ for the maps $\Phi_{q_1,q_2,r}$(a), $\Phi_{q_1,q_1,q_2,r}$(b), $\Phi_{q_1,q_2,q_2,r}$(c) and $\Phi_{q_1,q_2,q_1,r}$(d). We make the estimations using a orbit of 10.000 points. }
\label{top2}
\end{center}
\end{figure}

Following the notation of Section \ref{several}, we compute the topological entropy of the maps $\Phi_{q_1,q_2,r}$, $\Phi_{q_1,q_1,q_2,r}$, $\Phi_{q_1,q_2,q_2,r}$ and $\Phi_{q_1,q_2,q_1,r}$. We follow the same strategy used in the stability section, fixing $r=3.56$ and computing both the topological entropy and the Lyapunov exponents of these maps. We take a slightly bigger value of $r$ because the region of positive entropy for $3.5$ is really small. Note that for $h(f_{3.6})=0$. We also present the results for $(q_1,q_2)\in (0,4]^2$ because for bigger parameter regions the set of parameter values with positive topological entropy are located in a narrow region close to the axis. The results are shown in Figures \ref{top1} and \ref{top2}.

Figures \ref{top1} and \ref{top2} show again the existence of Parrondo's paradox. But here, it is important to realize that the figures are not symmetric with respect to the diagonal line. This fact indicates that the complexity depends on the order we iterate the $q$-deformations. This fact is different from the results from \cite{29} in which the complexity was independent of the order in which the iterations of the $q$-deformations were made.

To finish, we should remark that, although the application of several $q$-deformations seems to decrease the complexity of the corresponding $q$-deformed logistic map, this process is not capable of eliminating chaos for all the parameters $r$. For instance, for $r=4$, we have that $\Phi_{q_{1},\dots,q_{k},4}([0,\frac{1}{2}])=\Phi_{q_{1},\dots, q_{k},4}([\frac{1}{2},1])=I$, and so $\Phi_{q_{1},\dots,q_{k},4}$ has a $2$-horseshoe. Thus, the topological entropy of $\Phi_{q_{1},\dots,q_{k},4}$ is $\log(2)$ (see \cite[chapter 4]{30}). As the topological entropy is continuous for unimodal maps with positive entropy, there must be a neighborhood of $\Phi_{q_{1},\dots,q_{k},r}$ in which the maps $\Phi_{q_{1},\dots,q_{k},r}$ display positive entropy as well. So, this neighborhood becomes smaller when $k$ increases.

\section{Conclusions and future work}

In this paper, we extended the work of \cite{gupta} and, following the outline of \cite{29}, analyzed the $q$-deformed logistic map under a $q$-deformation, which displays worse properties than that of \cite{29}. As a result, we obtain that when several $q$-deformations are applied, the order in which we iterate them plays a role in bringing chaotic dynamics. Moreover, that order also influences the stability regions of the fixed points and, in general, the whole dynamics of the $q$-deformed family. The results of our study motivate us to raise two questions as follows. The first one is to make a comparative study of the $q$-deformed logistic family under other possible $q$-deformations (see, e.g. \cite{10,gupta}) and in a more general framework when order-preserving homeomorphisms of the interval $I=[0,1]$ are considered. The second problem is to explore other $q$-deformed models, i.e., replacing the logistic family for other possible ones.

\section*{Acknowledgements}
Jose S. C\'anovas has been supported by the grant MTM2017-84079-P funded by MCIN/AEI/10.13039/501100011033 and by "ERDF A way of making Europe", by the European Union.

Houssem Eddine Rezgui was supported by the research unit: ``Dynamical systems
and their applications", (UR17ES21), Ministry of Higher Education and
Scientific Research, Faculty of Science of Bizerte, Bizerte, Tunisia. This
work was made during the Houssem Eddine Rezgui's visit at Universidad Polit\'ecnica
de Cartagena under the Erasmus program KA107. The support of this
university is also gratefully acknowledged.

\end{document}